# Numerical Simulation of Two Dimentional sine-Gordon Solitons Using the Modified Cubic B-Spline Differential Quadrature Method

*H. S. Shukla*[1], *Mohammad Tamsir*[1*], *Vineet K. Srivastava*[2]

[1]Department of Mathematics & Statistics, DDU Gorakhpur University, Gorakhpur-273009, India
[2]ISRO Telemetry, Tracking and Command Network (ISTRAC), Bangalore-560058, India

**ABSTRACT**
In this article, a numerical simulation of two dimensional nonlinear sine-Gordon equation with Neumann boundary condition is obtained by using a composite scheme referred to as a modified cubic B-spline differential quadrature method (MCB-DQM). The modified cubic B-spline serves as a basis function in the differential quadrature method to compute the weighting coefficients. Thus, the sine-Gordon equation is converted into a system of second order ordinary differential equations (ODEs). We solve the resulting system of ODEs by an optimal five stage and fourth-order strong stability preserving Runge–Kutta scheme (SSP-RK54). Both damped and undamped cases are considered for the numerical simulation with Josephson current density function $\phi(x,y)=-1$ . The computed results are found to be in good agreement with the exact solutions and other numerical results available in literature.

## 1. Introduction

Consider the two-dimensional nonlinear sine-Gordon equation:

$$\frac{\partial^2 u}{\partial t^2}+\beta\frac{\partial u}{\partial t}=\frac{\partial^2 u}{\partial x^2}+\frac{\partial^2 u}{\partial y^2}-\phi(x,y)\sin(u),\quad x\in R,\quad t\ge 0, \tag{1.1}$$

where $R=\{(x,y): a\le x\le b\ ,\ \mathrm{c}\le y\le d\}$ , with the suitable initial conditions:

$$u(x,y,0)=f_1(x,y),\quad (x,y)\in R,$$
$$\frac{\partial u(x,y,0)}{\partial t}=f_2(x,y),\ (x,y)\in R, \tag{1.2}$$

and Neumann boundary conditions:

$$\frac{\partial u}{\partial n}=g_1(x,y,t),\ \text{ for } x=a \text{ and } x=b\ ,\ \mathrm{c}\le \mathrm{y}\le d,\ \mathrm{t}\ge 0$$
$$\frac{\partial u}{\partial n}=g_2(x,y,t),\ \text{ for } \mathrm{y}=c \text{ and } \mathrm{y}=d\ ,\ \mathrm{a}\le x\le b,\ \mathrm{t}\ge 0. \tag{1.3}$$

The parameter $\beta$ is known as dissipative term, which is supposed to be a real number with $\beta \geq 0$ . Eq. (1.1) reduces to the two-dimensional undamped sine-Gordon equation for $\beta = 0$, and to a damped one when $\beta > 0$. The function $\phi(x, y)$ represents as the Josephson current density, while $f_1(x, y)$ and $f_2(x, y)$ characterize wave modes or kinks and velocity, respectively. The kinks operates across the lines or closed curves in the $xy$-plane. These curves are known as the lines and ring soliton waves, respectively. Soliton waves appear in various applications and physical phenomena such as tsunamis, optical fiber signals, superconductivity, atmospheric waves and gravitational fields with cylindrical symmetry and so on. The simulation of the soliton waves represented by the sine-Gordon equation are of the great importance in research, and the numerical methods that are suitable for the success of the simulation of waves has led to a constant demand for achieving the higher accuracy. The Sine–Gordon equation appears in the propagation of fluxons in Josephson junctions, dislocations in crystals, solid state physics, nonlinear optics, stability of fluid motions and the motion of a rigid pendulum attached to a stretched wire, see [1-4]. Josephson junction model [4] consists of two superconducting layers separated by isolated barriers. This model can be described by the two-dimensional undamped sine-Gordon equation and has various applications in physics, electronics etc. Djidjeli et al. [7] modeled a soliton-like structure in higher dimensions.

In recent years, various numerical approaches have been developed for the solution of the sine-Gordon equation. First of all, Guo et al. [6] proposed two different difference schemes, namely, explicit and implicit schemes. Christiansen and Lomdahl [8] used a generalized leapfrog method for the numerical computation of two-dimensional undamped sine-Gordon equation while a finite elements method was used by Argyris et al. [9]. Both methods were successfully applied using the appropriate initial conditions with the latter one gave slightly more accurate results. Xin [10] modeled light bullets with two-dimensional sine–Gordon equation. Light bullets contain only a few electromagnetic oscillations under their envelopes and propagate long distances without essentially changing shapes. Author of [10] performed a modulation analysis and obseved that the sine-Gordon pulse envelopes undergo focusing–defocusing cycles. The evolution of lump and ring solitions of the two dimensional sine-Gordon equation and the evolution of standing and travelling breather-type waves are studied in Minzoni et al. [11,12]. Morever, lump and ring solitions can be applied to the Baby Skyrme model and to study the vortex models [11]. Sheng et al. [13] used a split cosine approach for the numerical simulation of two-dimensional sine-Gordon solitons.

Bratsos [15] proposed a modified predictor-corrector scheme. In [16], he used the method of lines while a three-time level fourth-order explicit finite-difference scheme used in [17] by him for solving two-dimensional sine-Gordon equation. The resulting nonlinear system was solved by a predictor–corrector (P-C) scheme using the explicit method as predictor and the implicit as corrector. An explicit numerical scheme and an improved numerical scheme for the numerical solutions of (2 +1) dimension sine-Gordon equation is proposed in [14,18]. Dehghan and his colleagues [19-24] proposed number of schemes for the numerical solutions of two-dimensional damped and undamped sine-Gordon equations. In [19], the authors proposed the dual reciprocity boundary element method and studied the various cases of line and ring solitons of the circular and elliptic shapes. In [20], the author proposed a numerical scheme based on the collocation points and thin plate splines radial basis functions. In [21], the author employed a continuous linear element approximation to study the boundary element solution. The non-linearity was dealt by a dual reciprocity and a predictor–corrector scheme, and non-linear and in-homogenous terms were converted into the boundary. In [22-24], the authors proposed various meshless schemes for the numerical solutions of the two-dimensional damped and undamped sine-Gordon equations. In [22], the author presented an implementation of the meshless local boundary integral equation (LBIE) scheme. In [23], proposed a local weak meshless scheme based on the radial point interpolation method (RPIM), while in [24], developed a local Petrov–Galerkin (MLPG) method by dealing the non-linearity terms by a predictor-corrector scheme while time-stepping scheme was used for the time derivative. Jiwari et al. [26] proposed a polynomial differential quadrature scheme to solve damped and undamped two-dimensional sine-Gordon equations and discussed various cases involving line solitons and ring solitons.

Differential quadrature technique can be seen in [26-32] for the solution of the linear and nonlinear one and two dimensional differential equations, while the cubic B-spline and sinc differential quadrature methods can be referred in [33, 34]. Arora & Singh [35] proposed a new method, namely, a "modified cubic B-spline differential quadrature method (MCB-DQM) to solve the one dimensional Burgers' equation and checked its efficiency and accuracy, and observed that the implementation of the method is very easy, powerful and efficient as compared to other existing numerical methods for the Burgers' equation. Recently, Jiwari & Yuan [36] extended MCB-DQM to show the computational modeling of two-dimensional reaction–diffusion Brusselator model with appropriate initial and Neumann boundary conditions. Authors of [37-39] developed an optimal strong stability preserving (SSP) high order time discretization schemes for the numerical solution of partial differential equations. SSP methods preserve the strong stability properties in any norm, semi norm or convex functional of the spatial discretization coupled with the first order Euler time stepping. They also discussed the description of the optimal explicit and implicit SSP Runge-Kutta and multistep methods.

In this work, we study the numerical simulation of the two-dimensional damped and undamped sine-Gordon equation for various cases including line and ring solitons of the circular and elliptic shapes using the MCB-DQM scheme. The efficacy and accuracy of the method is confirmed by taking test problem 4.1 of the undamped sine-Gordon equation with Neumann boundary conditions having the exact solution, which shows that the MCB-DQM results are acceptable and in good agreement with earlier studies available in literature. Test problems 4.2 to 4.6 shows the numerical simulation of damped and undamped sine-Gordon equation. Obtained numerical results are compared well with the numerical results available in litature and a good agreement is found. Rest of the article is organized as follows: In Section 2, the modified cubic B-spline differential quadrature method is introduced. In Section 3, the implementation procedure for the problem (1.1) – (1.3) is illustrated; in Section 4, various numerical examples are given to establish the applicability and accuracy of the method, while the Section 5 concludes our study.

## 2. Modified cubic B-spline differential quadrature method in two dimension

The differential quadrature method (DQM) was firstly introduced by Bellman et al. [26] in 1972. This method approximates the spatial derivatives of a function using the weighted sum of the functional values at the certain discrete points. In DQM, the weighting coefficients are computed by using several kinds of test functions such as spline function [33], sinc function [34], Lagrange interpolation polynomials and Legendre polynomials [27-32] etc. This section re-describes the MCB-DQM [35,36] to complete our problem in two dimension. It is assumed that the $M$ and $N$ grid points: $a = x_1 < x_2, ..... < x_M = b$ and $c = y_1 < y_2, ..... < y_N = d$ are distributed uniformly with the step size $\Delta x = x_{i+1} - x_i$ and $\Delta y = y_{j+1} - y_j$ in $x$ and $y$ directions, respectively. According to the method, the first and second order spatial partial derivatives of $u(x, y, t)$ with respect to $x$ (keeping $y_j$ as fixed), approximated at $x_i$, are defined as:

$$\frac{\partial u(x_i, y_j, t)}{\partial x} = \sum_{k=1}^{N} w_{ik}^{(1)} u(x_k, y_j, t), \qquad i = 1, 2, ..., M \tag{2.1}$$

$$\frac{\partial^2 u(x_i, y_j, t)}{\partial x^2} = \sum_{k=1}^{N} w_{ik}^{(2)} u(x_k, y_j, t), \qquad i = 1, 2, ..., M. \tag{2.2}$$

Similarly, the first and second order spatial partial derivatives of $u(x, y, t)$ with respect to $y$ (keeping $x_i$ as fixed), approximated at $y_j$ , are defined as:

$$\frac{\partial u(x_i, y_j, t)}{\partial y} = \sum_{k=1}^{M} \bar{w}_{jk}^{(1)} u(x_i, y_k, t), \qquad j = 1, 2, ..., N \tag{2.3}$$

$$\frac{\partial^2 u(x_i, y_j, t)}{\partial y^2} = \sum_{k=1}^{M} \bar{w}_{jk}^{(2)} u(x_i, y_k, t), \qquad j = 1, 2, ..., N, \tag{2.4}$$

where $w_{ij}^{(r)}$ and $\bar{w}_{ij}^{(r)}$ , $r = 1, 2$ are the weighting coefficients of the $r$th-order spatial partial derivatives with respect to $x$ and $y$.

The cubic B-spline basis functions [32] at the knots are defined as:

$$\varphi_m(x) = \frac{1}{h^3} \begin{cases} (x - x_{m-2})^3, & x \in (x_{m-2}, x_{m-1}) \\ (x - x_{m-2})^3 - 4(x - x_{m-1})^3, & x \in (x_{m-1}, x_m) \\ (x_{m+2} - x)^3, & x \in (x_{m+1}, x_{m+2}) \\ 0, & \text{otherwise}, \end{cases} \tag{2.5}$$

where $\{\varphi_0, \varphi_1, ..., \varphi_N, \varphi_{N+1}\}$ is chosen in such a way that it forms a basis over the domain $R = \{(x, y) : a \le x \le b;\ c \le y \le d\}$. The values of cubic B-splines and its derivatives at the nodal points are depicted in **Table 1**.

**Table 1:** Coefficients of the cubic B-spline $\varphi_m$ and its derivatives at the node $x_m$.

| | $x_{m-2}$ | $x_{m-1}$ | $x_m$ | $x_{m+1}$ | $x_{m+2}$ |
|---|---|---|---|---|---|
| $\varphi_m(x)$ | 0 | 1 | 4 | 1 | 0 |
| $\varphi'_m(x)$ | 0 | $3/h$ | 0 | $-3/h$ | 0 |
| $\varphi''_m(x)$ | 0 | $6/h^2$ | $-12/h^2$ | $6/h^2$ | 0 |

Now, we modify the cubic B-spline basis functions in such a way that the resulting matrix system of the equations becomes diagonally dominant. Then, the modified cubic B-splines [35] are defined as:

$$\left.\begin{aligned} &\phi_1(x) = \varphi_1(x) + 2\varphi_0(x) \\ &\phi_2(x) = \varphi_2(x) - \varphi_0(x) \\ &\phi_m(x) = \varphi_m(x) \; for \; \mathrm{m} = 3, ..., N - 2 \\ &\phi_{N-1}(x) = \varphi_{N-1}(x) - \varphi_{N+1}(x) \\ &\phi_N(x) = \varphi_N(x) + 2\varphi_{N+1}(x) \end{aligned}\right\}, \tag{2.6}$$

where $\{\phi_1, \phi_2, ..., \phi_N\}$ forms a basis over the domain $R = \{(x, y) : a \le x \le b; \; c \le y \le d\}$. Since, in Eq. (2.1), the $x$-axis is fixed, therefore, by substituting the values of $\phi_m(x)$, $m = 1, 2, ..., N$ in Eq. (2.1), we find the weighing coefficient $w_{ik}^{(1)}$ by the equation:

$$\phi_m'(x_i) = \sum_{k=1}^{N} w_{ik}^{(1)} \phi_m(x_k), \; \mathrm{i} = 1, 2, ..., M. \tag{2.7}$$

From Eq. (2.5), (2.6) and Table 1, Eq. (2.7) reduces into a tridiagonal system of equations:

$$A\vec{w}^{(1)}[i] = \vec{R}[i], \; for \; i = 1, 2, ..., M, \tag{2.8}$$

where $\vec{w}^{(1)}[i] = \left[ w_{i1}^{(1)}, w_{i2}^{(1)}, ..., w_{iN}^{(1)} \right]^T$ is the weighting coefficient vector corresponding to $x_i$, $\vec{R}[i] = \left[ \phi_{1,i}', \; \phi_{2,i}', \; ......, \; \phi_{N-1,i}', \; \phi_{N,i}' \right]^T$, and the coefficient matrix $A$ is given by:

$$A = \begin{bmatrix} \phi_{1,1} & \phi_{1,2} & & & & & \\ \phi_{2,1} & \phi_{2,2} & \phi_{2,3} & & & & \\ & \phi_{3,2} & \phi_{3,3} & \phi_{3,4} & & & \\ & & \ddots & \ddots & \ddots & & \\ & & & \phi_{N-2,N-3} & \phi_{N-2,N-2} & \phi_{N-2,N-1} & \\ & & & & \phi_{N-1,N-2} & \phi_{N-1,N-1} & \phi_{N-1,N} \\ & & & & & \phi_{N,N-1} & \phi_{N,N} \end{bmatrix} \tag{2.9}$$

We notice that the coefficient matrix $A$, being the diagonally dominant, is invertible. The tridiagonal system of linear equations (2.8) is solved for each $i$ using the Thomas algorithm. This gives the weighting coefficients $w_{ik}^{(1)}$ of the first order partial derivative. In this way, we obtain the weighting coefficients $w_{ij}^{(2)}, 1 \le i, j \le N$. The weighting coefficients $w_{ij}^{(2)}, 1 \le i, j \le N$, are determined by using the formula [32]:

$$\begin{cases} w_{ij}^{(r)} = r\left( w_{ij}^{(1)} w_{ii}^{(r-1)} - \dfrac{w_{ij}^{(r-1)}}{x_i - x_j} \right), \; for \;\; i \neq j \text{ and } i = 1,2,3,...,N; \quad r = 2,3,...,N-1 \\ w_{ii}^{(r)} = -\sum\limits_{j=1,\, j\neq i}^{N} w_{ij}^{(r)}, \;\; \text{for } i = j, \end{cases} \tag{2.10}$$

where $\bar{w}_{ij}^{(r-1)}$ and $\bar{w}_{ij}^{(r)}$ are the weighting coefficients of the $(r-1)^{\text{th}}$ and $r^{\text{th}}$ order partial derivatives with respect to $x$.

In the same way, we obtain weighting coefficients $\bar{w}_{jk}^{(1)}$ of the first order partial derivatives with respect to $y$ by using the modified cubic B-Spline functions in Eq. (2.3). The weighting coefficients $\bar{w}_{ij}^{(2)}, 1 \leq i, j \leq N$ for the second derivatives can be calculated by the recurrence formule:

$$\begin{cases} \bar{w}_{ij}^{(r)} = r\left( \bar{w}_{ij}^{(1)} \bar{w}_{ii}^{(r-1)} - \dfrac{\bar{w}_{ij}^{(r-1)}}{x_i - x_j} \right), \; for \;\; i \neq j \text{ and } i = 1,2,3,...,N; \quad r = 2,3,...,N-1 \\ \bar{w}_{ii}^{(r)} = -\sum\limits_{j=1,\, j\neq i}^{N} \bar{w}_{ij}^{(r)}, \;\; \text{for } i = j, \end{cases} \tag{2.11}$$

where $\bar{w}_{ij}^{(r-1)}$ and $\bar{w}_{ij}^{(r)}$ are the weighting coefficients of the $(r-1)^{\text{th}}$ and $r^{\text{th}}$ order partial derivatives with respect to $y$.

## 3. Implementation of Method to two-dimensional sine-Gordon equation

On substituting the approximate values of the spatial derivatives computed by the MCB-DQM, Eq. (1.1) can be re-written as:

$$\frac{\partial^2 u(x_i, y_j, t)}{\partial t^2} + \beta \frac{\partial u(x_i, y_j, t)}{\partial t} = \sum_{k=1}^{M} w_{ik}^{(2)} u(x_k, y_j) + \sum_{k=1}^{N} \bar{w}_{jk}^{(2)} u(x_i, y_k) - \phi(x_i, y_j) \sin\left(u(x_i, y_j)\right),$$

$$(x_i, y_j) \in R, \; t > 0, \; i = 1,2,...,M, \quad j = 1,2,...,N, \tag{3.1}$$

with the initial conditions:

$$\begin{aligned} & u(x_i, y_j, 0) = f_1(x_i, y_j) \\ & \frac{\partial u(x_i, y_j, 0)}{\partial t} = f_2(x_i, y_j). \end{aligned} \tag{3.2}$$

Eq. (3.1) reduces into a system of second order differential equations:

$$\frac{d^2u\left(x_i,y_j,t\right)}{dt^2}=L\left(u\left(x_i,y_j,t\right)\right),\ \ i=1,2,...,M \text{ and } j=1,2,...,N, \tag{3.3}$$

where $L$ is a nonlinear differential operator.

The approximation of Neumann boundary conditions (1.3) at $x=a$ and $x=b$ is given as [25,36]:

$$\sum_{k=1}^{M} w_{1,k}^{(1)} u_{k,j} = g_1 \tag{3.4}$$

$$\sum_{k=1}^{M} w_{M,k}^{(1)} u_{k,j} = g_2 \tag{3.5}$$

On solving Eqs. (3.4) and (3.5) for $u_{1,j}$ and $u_{M,j}$, we have

$$u_{1,j}=\frac{w_{1,M}^{(1)}\left(g_2-S_2\right)-w_{M,M}^{(1)}\left(g_1-S_1\right)}{\left(w_{1,M}^{(1)}w_{M,1}^{(1)}-w_{1,1}^{(1)}w_{M,M}^{(1)}\right)}, \tag{3.6}$$

$$u_{M,j}=\frac{w_{M,1}^{(1)}\left(g_1-S_1\right)-w_{1,1}^{(1)}\left(g_2-S_2\right)}{\left(w_{1,M}^{(1)}w_{M,1}^{(1)}-w_{1,1}^{(1)}w_{M,M}^{(1)}\right)}, \tag{3.7}$$

where $S_1=\sum_{k=2}^{M-2} w_{1,k}^{(1)}u_{k,j},\ S_2=\sum_{k=2}^{M-2} w_{M,k}^{(1)}u_{k,j},\ g_1=g_1\left(a,y,t\right),\ g_2=g_2\left(b,y,t\right)$ and $j=1,\ 2,....,\ N.$

Similary, the approximation of the Neumann boundary conditions (1.3) at $y=c$ and $y=d$ is given by:

$$\sum_{k=1}^{N} \bar{w}_{1,k}^{(1)} u_{i,k} = g_3\left(x,c,t\right), \tag{3.8}$$

$$\sum_{k=1}^{N} \bar{w}_{N,k}^{(1)} u_{i,k} = g_4\left(x,d,t\right). \tag{3.9}$$

On solving Eqs. (3.8) and (3.9) for $u_{i,1}$ and $u_{i,N}$, we have

$$u_{i,1}=\frac{\bar{w}_{1,N}^{(1)}\left(g_4-S_4\right)-\bar{w}_{N,N}^{(1)}\left(g_3-S_3\right)}{\left(\bar{w}_{1,N}^{(1)}w_{N,1}^{(1)}-\bar{w}_{1,1}^{(1)}\bar{w}_{N,N}^{(1)}\right)}, \tag{3.10}$$

$$u_{i,N}=\frac{\bar{w}_{N,1}^{(1)}\left(g_3-S_3\right)-\bar{w}_{1,1}^{(1)}\left(g_4-S_4\right)}{\left(\bar{w}_{1,N}^{(1)}\bar{w}_{N,1}^{(1)}-\bar{w}_{1,1}^{(1)}\bar{w}_{N,N}^{(1)}\right)}, \tag{3.11}$$

where $S_3=\sum_{k=2}^{N-2} \bar{w}_{1,k}^{(1)}u_{i,k},\ S_4=\sum_{k=2}^{N-2} \bar{w}_{N,k}^{(1)}u_{i,k},\ g_3=g_3\left(x,c,t\right),\ g_4=g_4\left(x,d,t\right)$ and $i=1,\ 2,....,$ M.

The system of second order ordinary differential equations (3.3) together with the initial conditions (3.2) and Neumann boundary conditions (3.6), (3.7), (3.10) (3.11) are solved by SSP-RK54 scheme.

## 4. Results and discussion

In section, we consider six problems of the line and ring solitons to provide the MCB-DQM numerical solutions of Eq. (1.1) with the appropriate initial conditions (1.2) and Neumann boundary conditions (1.3). The accuracy and efficiency of the scheme is measured by evaluating the $L_\infty$ and RMS error norms, defined as:

$$\left.\begin{aligned} &L_\infty = \max_{\substack{1\le i\le M\\ 1\le j\le N}} |u_{i,j}^{exact} - u_{i,j}^{numerical}| \\ &\text{and} \quad \text{RMS} = \frac{1}{M\times N}\sqrt{\sum_{i=1}^{M}\sum_{j=1}^{N}\left|u_{i,j}^{exact} - u_{i,j}^{numerical}\right|^2} \end{aligned}\right\}. \tag{4.1}$$

Also, the computed results are compared with those available in published papers.

### 4.1. Test Problem

In this problem, numerical solution of the two-dimensional sine-Gordon equation (1.1) are obtained for the case $\phi(x,y) = -1$ and $\beta = 0$ in the region $-7 \le x, y \le 7$ with the initial conditions:

$$\begin{aligned} &u(x,y,0) = 4\tan^{-1}\left(\exp(x+y)\right), \qquad -7 \le x, y \le 7, \\ &u_t(x,y,0) = -\frac{4\exp(x+y)}{1+\exp(2x+2y)}, \qquad -7 \le x, y \le 7, \end{aligned} \tag{4.2}$$

and Neumann boundary conditions:

$$\begin{aligned} &\frac{\partial u}{\partial x} = -\frac{4\exp(x+y+t)}{\exp(2t)+\exp(2x+2y)}, \qquad \text{for } x = \text{-}7 \text{ and } x = 7,\ -7 \le y \le 7,\ t > 0, \\ &\frac{\partial u}{\partial y} = -\frac{4\exp(x+y+t)}{\exp(2t)+\exp(2x+2y)}, \qquad \text{for } y = \text{-}7 \text{ and } y = 7,\ -7 \le x \le 7,\ t > 0. \end{aligned} \tag{4.3}$$

The exact solution of Eq. (1.1) with the initial conditions (4.2) and for the case $\phi(x,y) = -1$ and $\beta = 0$ is given by:

$$u(x,y,t) = 4\tan^{-1}\left(\exp(x+y-t)\right).$$

The numerical solution of the test problem 4.1 are obtained in the region $-7 \le x, y \le 7$ using the parameters: the time step $\Delta t = 0.001$ and the space step size $\Delta x = \Delta y = 0.25$. **Table 2** shows the $L_\infty$ and root mean square (RMS) error norms at the different time levels $t = 1,\ 3,\ 5$ and $t = 7$. From **Table 2**, It can be seen that the MCB-DQM results are better than those of [7, 20, 25 ]. **Fig. 1** shows the comparison between the numerical and exact solutions at $t = 5$ and $t = 7$ with $\Delta t = 0.001$ and grid size $31 \times 31$ , which we notice that are resemble to each other.

**Table 2:** Comparison of the $L_\infty$ and RMS errors of problem 4.1 at the different time levels $t$.

| $t$ | $L_\infty$ − error | | | | RMS-error | | |
|---|---|---|---|---|---|---|---|
| | [7] | [20] | [25] | present | [20] | [25] | present |
| 1 | 0.0350 | 0.0670 | 0.0027 | 0.0003 | 0.0050 | 0.0005 | 0.0002 |
| 3 | 0.0431 | 0.0834 | 0.0020 | 0.0006 | 0.0103 | 0.00005 | 0.0004 |
| 5 | 0.0404 | 0.1015 | 0.0033 | 0.0008 | 0.0145 | 0.0007 | 0.0007 |
| 7 | 0.0353 | 0.1516 | 0.0059 | 0.0012 | 0.0187 | 0.0011 | 0.0010 |

**4.2.** Circular ring solitons

The circular ring solitons of the two-dimensional sine-Gordon equation (1.1) are obtained for the case $\phi(x, y) = -1$ and $\beta = 0$ in the region $-7 \le x, y \le 7$ with the initial conditions:

$$u(x,y,0) = 4\tan^{-1}\left(\exp\left(3 - \sqrt{x^2+y^2}\right)\right), \qquad -7 \le x, y \le 7,$$

$$u_t(x,y,0) = 0, \qquad -7 \le x, y \le 7, \tag{4.3}$$

and Neumann boundary conditions:

$$\left.\begin{aligned} \frac{\partial u}{\partial x} &= 0, \quad \text{for } x = \text{-}7 \text{ and } x = 7,\ -7 \le y \le 7,\ t > 0 \\ \frac{\partial u}{\partial y} &= 0, \quad \text{for } y = \text{-}7 \text{ and } y = 7,\ -7 \le x \le 7,\ t > 0 \end{aligned}\right\}. \tag{4.4}$$

Fig. 2. shows the initial conditions (i.e., $t = 0$) and the numerical solutions of the circular ring solitons at $t = 2.8,\ 5.6,\ 8.4,\ 11.2,$ and $t = 12.6$ in terms of $\sin(u/2)$. From Fig. 2 depicts that the ring soliton shrinks from initial position ($t = 0$) to $t = 2.8$ and appears as a single ring soliton. The expansion phase starts from $t = 5.6$ and contiues until $t = 11.2$, a radiation

appears at the beginning of the expansion phase. At $t = 11.2,$ a ring soliton is nearly reformed. Fig. 2 graphs are in good agreement with those given in published papers [7-9, 15, 18, 23, 25].

**4.3.** Elliptical ring solitons

The elliptical ring solitons of the two-dimensional sine-Gordon equation (1.1) are obtained for the case $\phi(x, y) = -1$ and $\beta = 0$ in the region $-7 \le x, y \le 7$ with the initial conditions:

$$u(x, y, 0) = 4\tan^{-1}\left(\exp\left(3 - \sqrt{\frac{(x-y)^2}{3} + \frac{(x+y)^2}{2}}\right)\right), \quad -7 \le x, y \le 7,$$

$$u_t(x, y, 0) = 0, \quad -7 \le x, y \le 7, \tag{4.5}$$

and Neumann boundary conditions (4.4).

Initial conditions (i.e., $t = 0$) and the numerical solutions of the elliptical ring solitons at $t = 1.6,\ 3.2,\ 4.8,\ 8.0,\ 9.6,\ 11.2$ in terms of $\sin(u/2)$ are shown in Fig. 3. From Fig. 3, we notice that the elliptical ring soliton shrinks from its initial position ($t = 0$) until $t = 3.2.$ The temporal behavior of this soliton wave consists of a shrinking and a blowing up phase. From $t = 8.0,$ the elliptical ring solitons begin to form a circular ring soliton and at $t = 11.2,$ almost a circular ring soliton is seen as [9, 16, 19, 21, 23, 25].

**4.4.** Elliptical breather

In this problem, the elliptical breathe for the case $\phi(x, y) = -1$ and $\beta = 0$ in the region $-7 \le x, y \le 7$ with the following initial conditions:

$$u(x, y, 0) = 4\tan^{-1}\left(2.0\operatorname{sec}h\left(0.866\sqrt{\frac{(x-y)^2}{3} + \frac{(x+y)^2}{2}}\right)\right), \quad -7 \le x, y \le 7,$$

$$u_t(x, y, 0) = 0, \quad -7 \le x, y \le 7, \tag{4.6}$$

and with Neumann boundary conditions (4.3) same as taken in circular ring solitons are discussed.

Fig. 4 depicts the initial conditions (i.e., $t = 0$) and the numerical solutions in three dimension ($3D$), and the contours plots for the elliptical breather at $t = 1.6,\ 8.0,\ 9.6,\ 11.2, 12.8,\ 14.8, 15.2$. These graphs reveal that the elliptical breather rotates clockwise from its initial position about the major axis $y = -x$ and shrinks until $t = 1.6$. At $t = 11.2$, the major axis is nearly recovered to its initial position but a tough oscillations are seen. At $t = 12.8$, an expansion phase begins. A comparison is made between our graphs and the graphs available in litature. A good agreement is found with those given in [8, 9, 16, 17, 23, 25].

**4.5.** Superposition of two orthogonal lines

In this problem, the superposition of two orthogonal line solitons are discussed for the case $\phi(x, y) = -1$ and $\beta = 0.05$ in the region $-6 \le x, y \le 6$ with the following initial conditions:

$$u(x, y, 0) = 4\tan^{-1}\left(\exp(x) + \exp(y)\right), \qquad -6 \le x, y \le 6,$$

$$u_t(x, y, 0) = 0, \qquad -7 \le x, y \le 7, \tag{4.6}$$

and with Neumann boundary conditions (4.4) .

The numerical solutions of this problem are obtained at the different time levels $t = 1,\ 3,\ 7,\ 10,\ 15,\ 20$ with $\Delta t = 0.001$ and for grid size $31\times 31$ and are shown in Fig. 5. It can be seen from the graphs that the break up of two orthogonal line solitons at $t = 1$ is found which moves away from each other and a separation occurs between them without any deformation until $t = 3$. At $t = 7$, a deformation is appeared. A good agreement is found with the results given in [25].

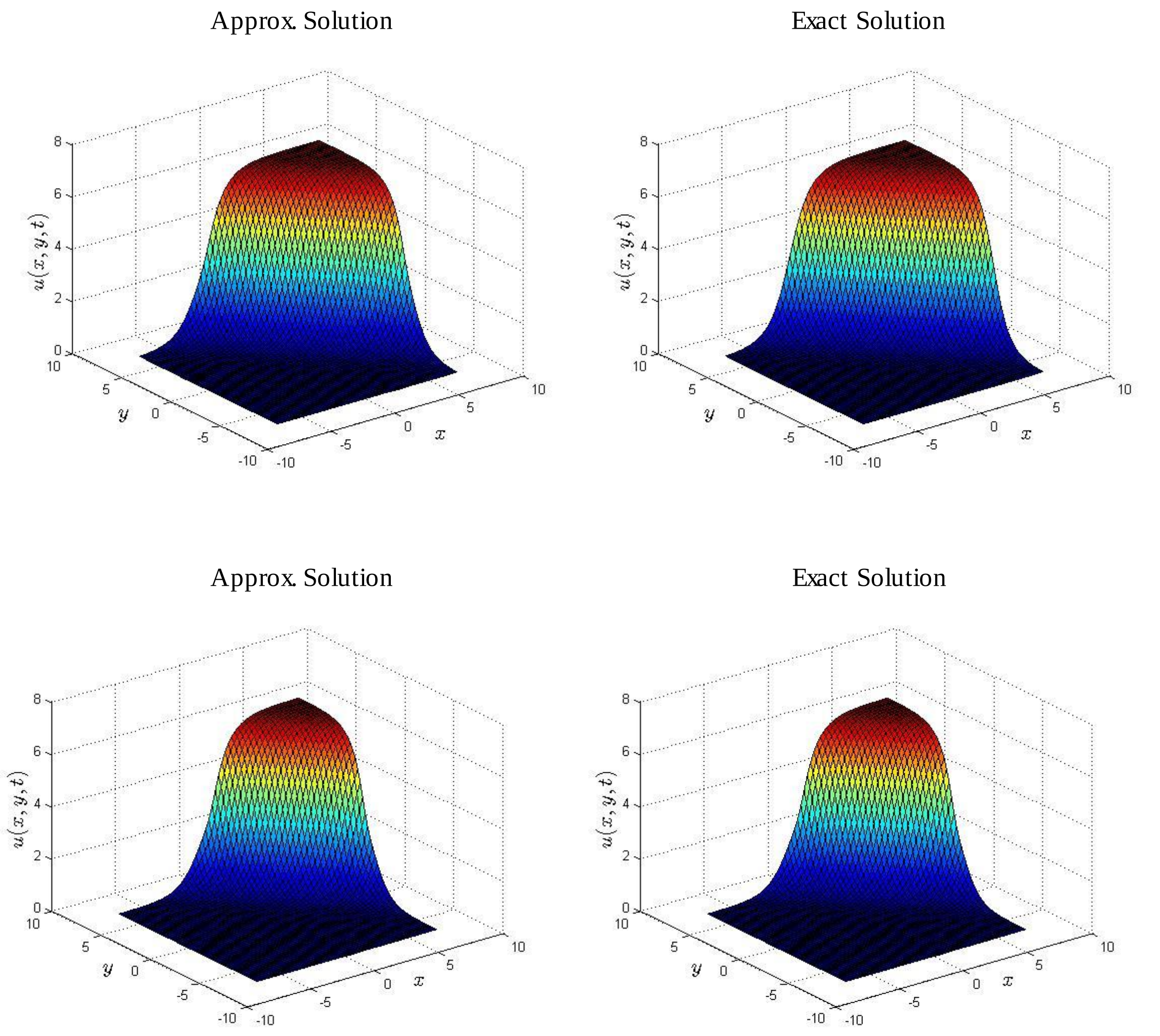

**Fig. 1.** Simulated and Exact solutions for $t = 5.0$ and $t = 7.0$ with $\Delta t = 0.001$, $\Delta x = \Delta y = 0.25$ for the problem 5.1.

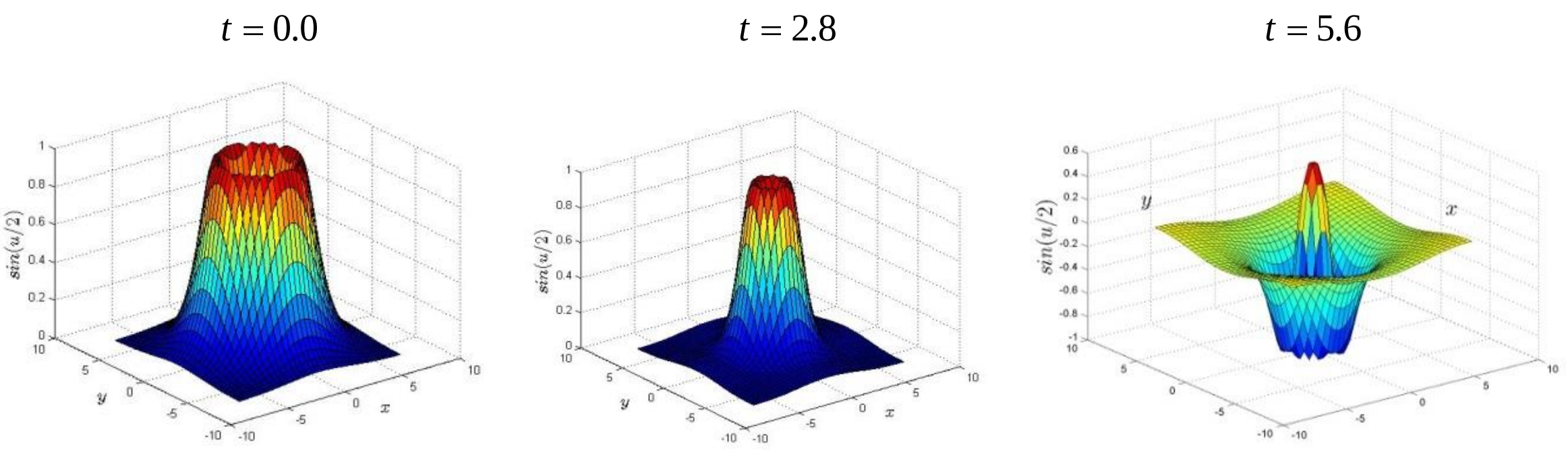

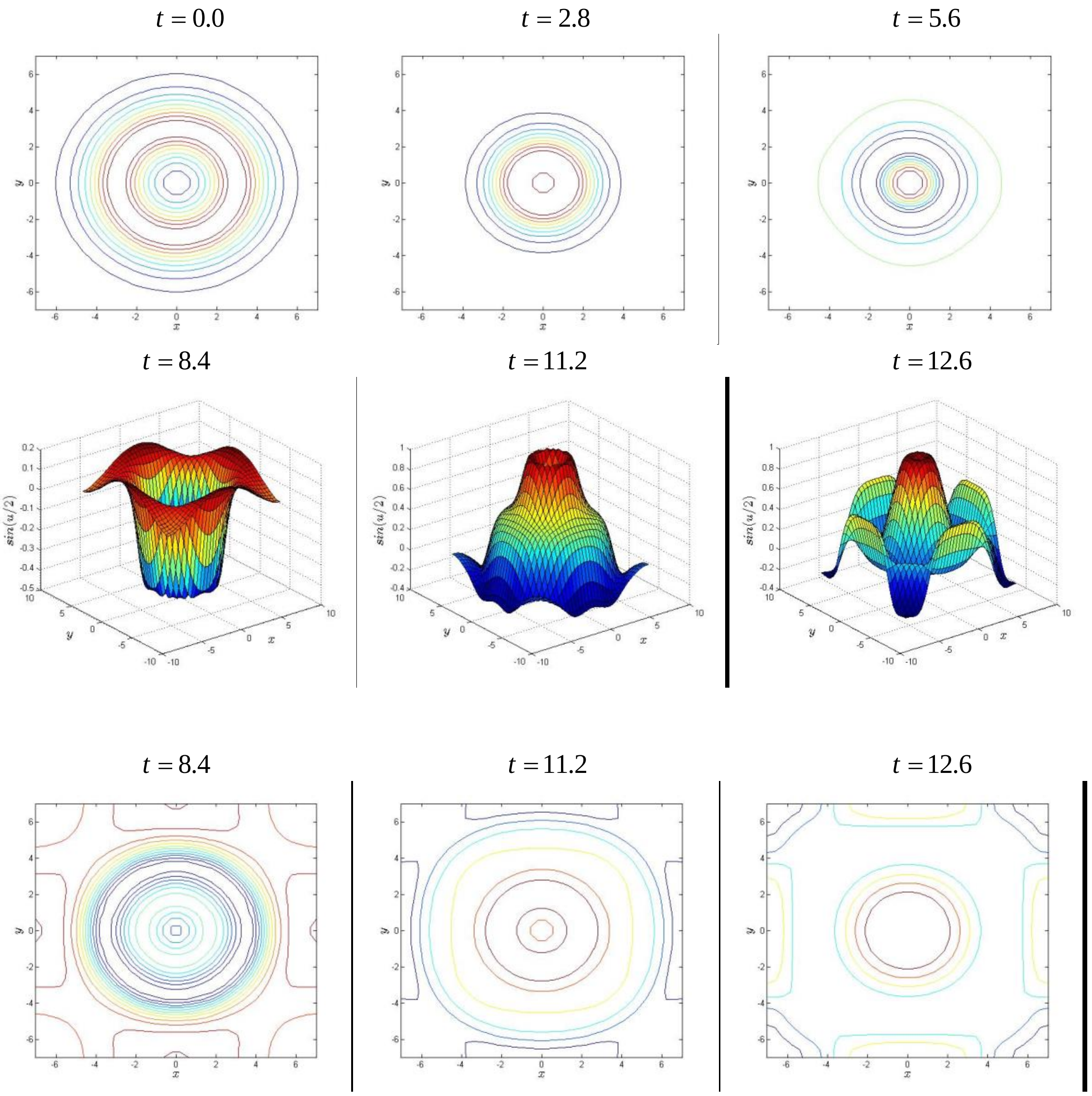

**Fig. 2.** Initial conditions and Numerical solutions of circular ring solitons at times $t = 2.8,\ 5.6,\ 8.4,\ 11.2,\ 12.6$ with $\Delta t = 0.2,\ \ \Delta x = \Delta y = 0.4$.

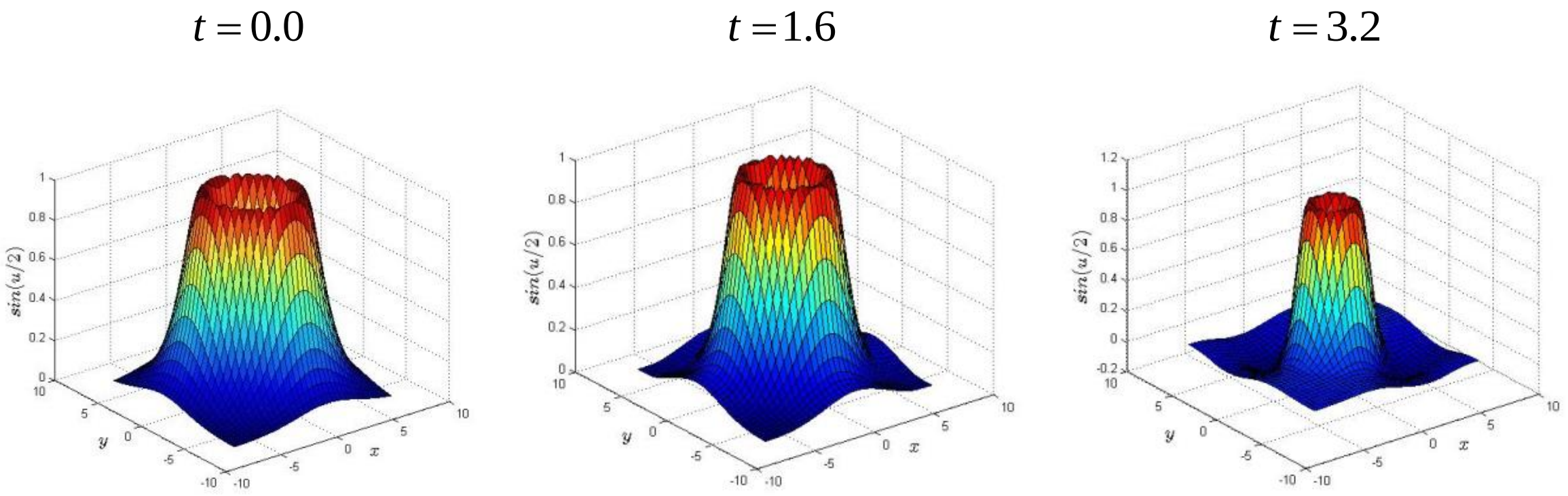

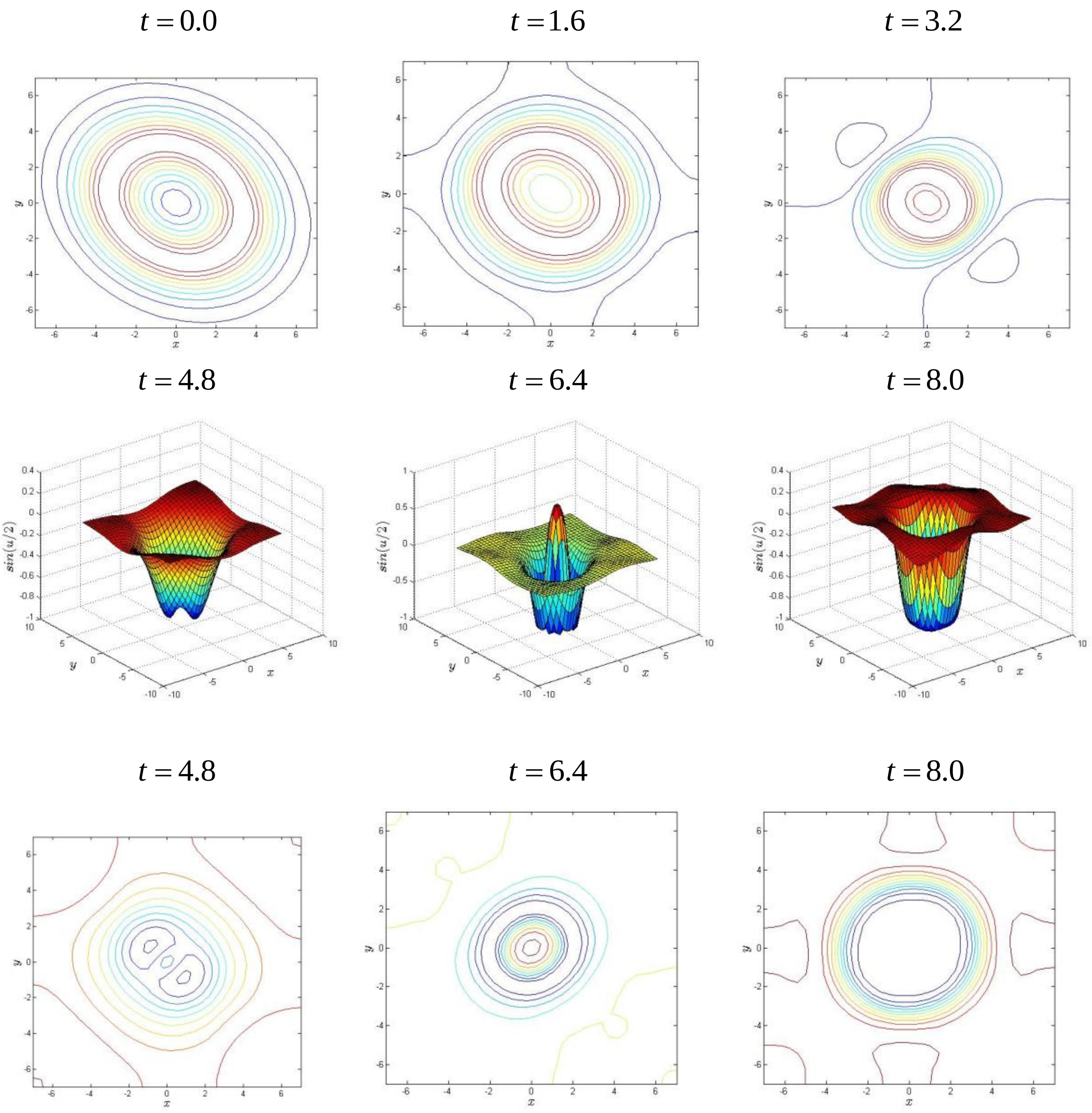

**Fig. 3.** Initial conditions and Numerical solutions of elliptical ring soliton at times $t = 1.6,\ 3.2,\ 4.8, 6.4,\ 8.0,\ 9.6, 11.2$ with $\Delta t = 0.2,\ \ \Delta x = \Delta y = 0.4$.

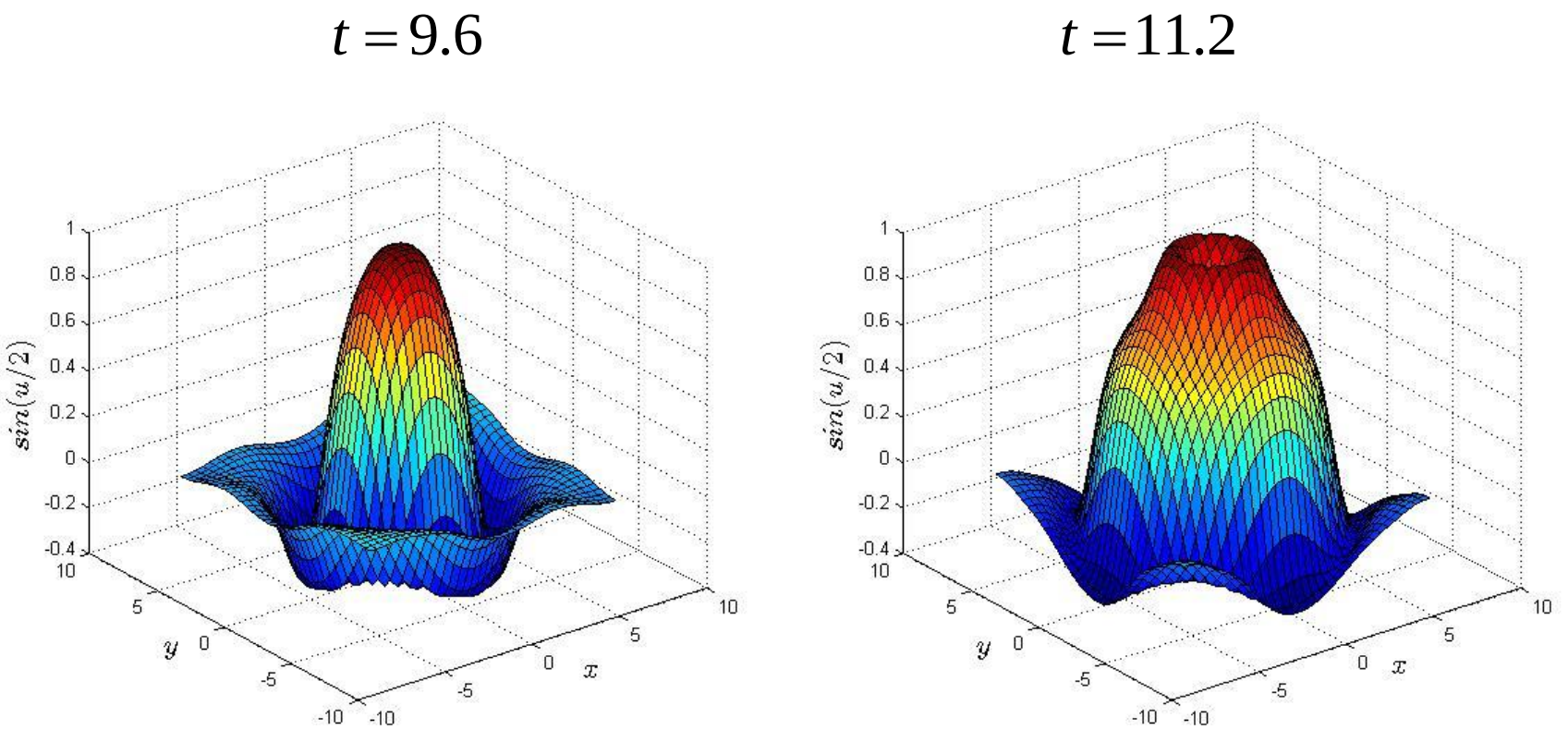

$t = 9.6$ $\qquad\qquad$ $t = 11.2$

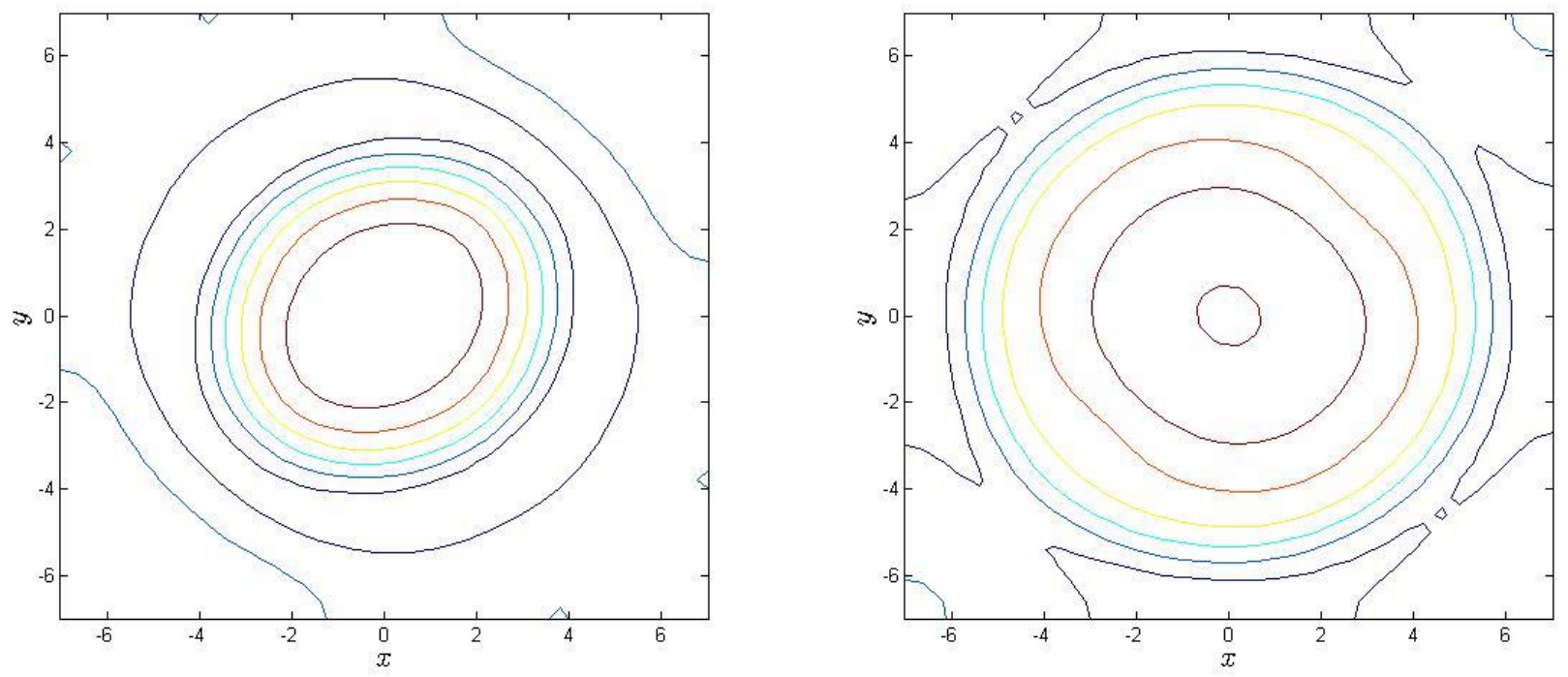

**Fig. 3.** Continued.

**4.6.** Line solitons in an inhomogeous medium

The numerical solutions for a line soliton in an inhomogeous medium are discussed for the case $\phi(x, y) = -1$ and $\beta = 0.05$ in the region $-7 \le x, y \le 7$ with the following initial conditions:

$$u(x, y, 0) = 4\tan^{-1}\left(\exp\left(\frac{x-3.5}{0.954}\right)\right), \qquad -7 \le x, y \le 7,$$

$$u_t(x, y, 0) = 0.629\sec h\left(\exp\left(\frac{x-3.5}{0.954}\right)\right), \qquad -7 \le x, y \le 7, \tag{4.7}$$

and with Neumann boundary conditions (4.4) .

The numerical solutions of this problem in an inhomogeous medium at the different time levels $t = 0,\ 6,\ 12,\ 18$ with $\Delta t = 0.001$ and for grid size $31\times 31$ are presented in Fig. 6. It is clear from the graphs that the line soliton moves slightly in a straight line during the transmission. A twist in its straightness appears as $t$ increases. Due to the inhomogeneity of the medium, this movement seems to be prevented. Finally, the soliton reforms its straightness at $t = 18$. Christiansen and Lomdahl [8] claimed that this effect is due to the boundary conditions. The graphs are compared with those given in [20, 25] and a good agreement is found.

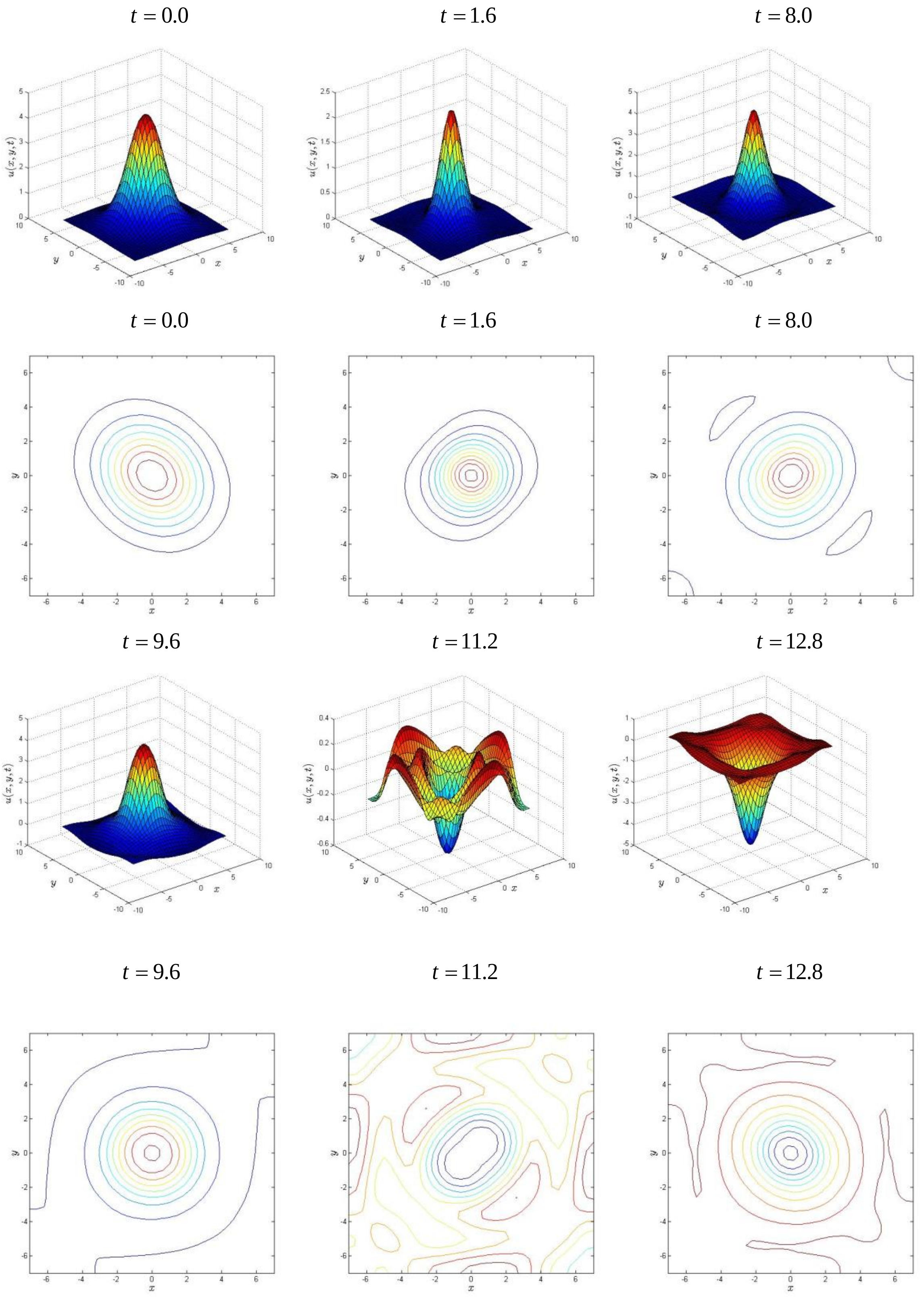

**Fig. 4.** Initial conditions and Numerical solutions of elliptical breather at times $t = 1.6,\ 8.0,\ 9.6,\ 11.2,\ 12.8,\ 14.8,\ 15.2$ with $\Delta t = 0.2,\ \ \Delta x = \Delta y = 0.4$.

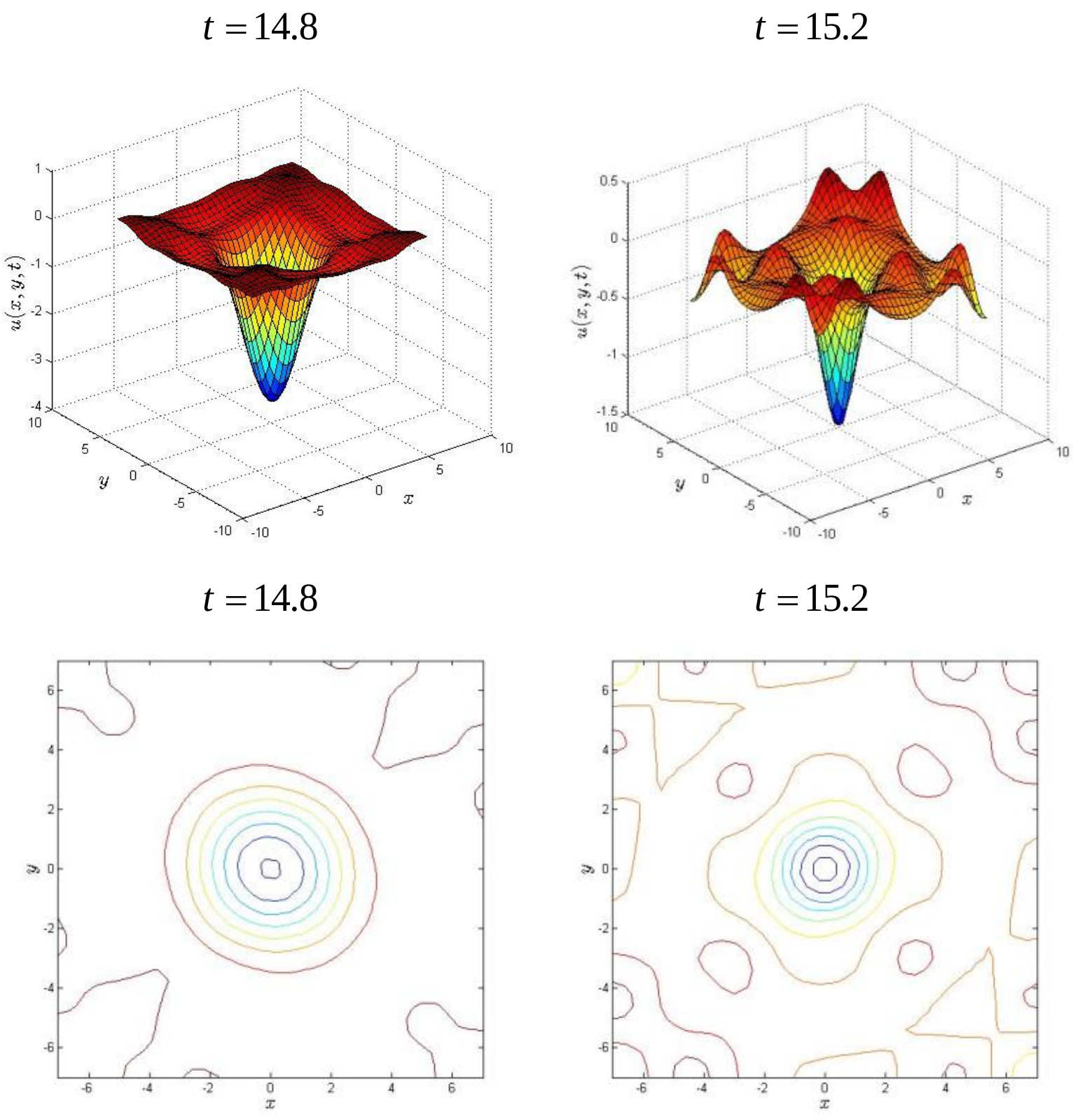

**Fig. 4.** Continued.

## 5. Conclusions

In this study, we discussed a composite scheme: MCB-DQM, in space with SSP-RK54 scheme, in time for solving the two-dimensional damped/undamped sine-Gordon equation numerically. In the study, various cases of the line and ring solitions are discussed together with circular and elliptical shapes. The computed numerical results and graphs are compared with those available in liturature. From the first test problem, it is evident that the described scheme provides more accurate results, and an excellent numerical approximation to the exact solution is achieved. For the rest of the problems. From the three dimensional and contour plots, it is demonstrated that the scheme gives analogous results as given in liturature. Additionally, we notice that the present scheme is very easy to apply and gives more accuracy even for small grid sizes.

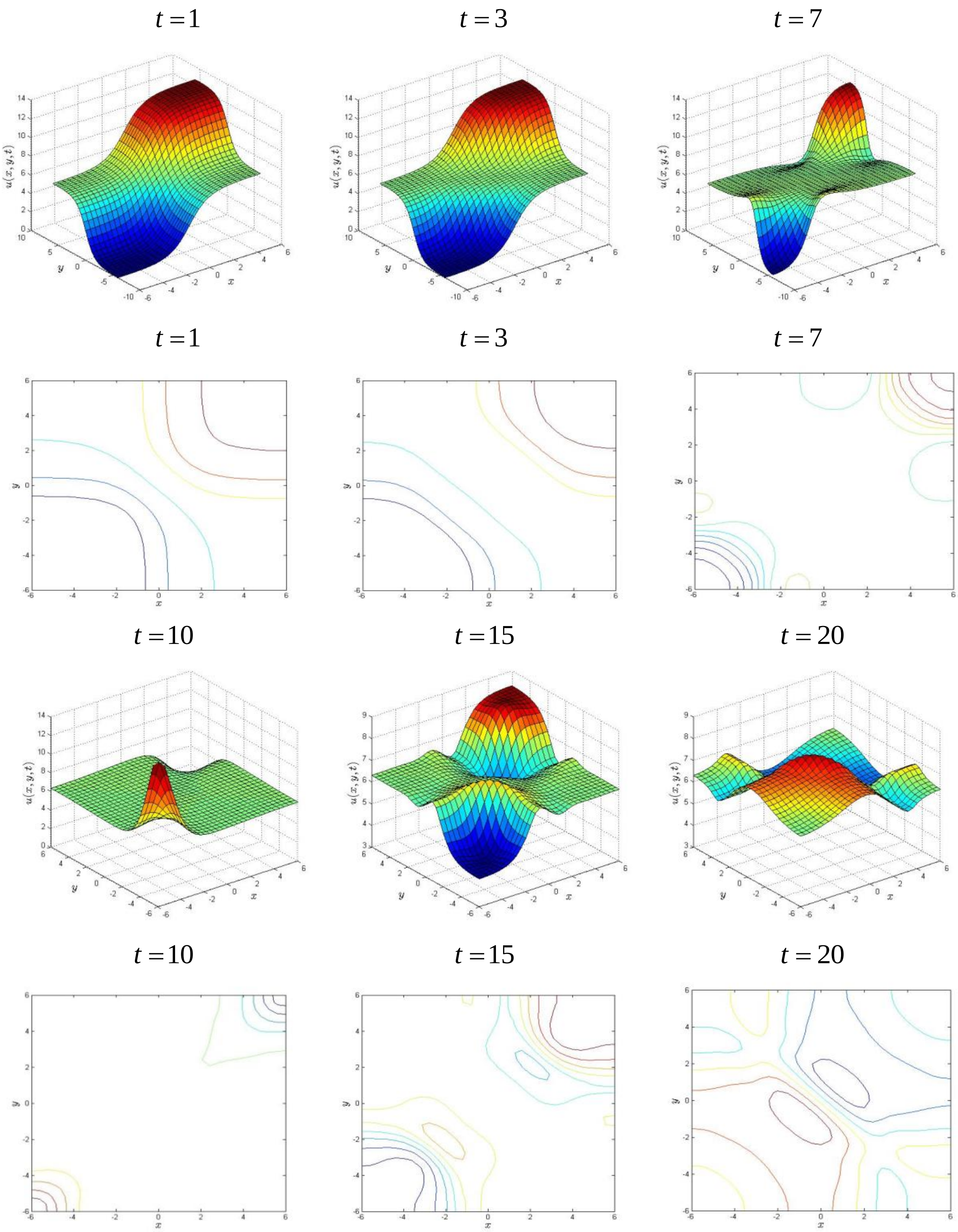

**Fig. 5.** Numerical solutions of superposition of two orthogonal lines at times $t = 1,\ 3,\ 7,\ 10,\ 15,\ 20$ with $\Delta t = 0.001,\ \Delta x = \Delta y = 0.4$ .

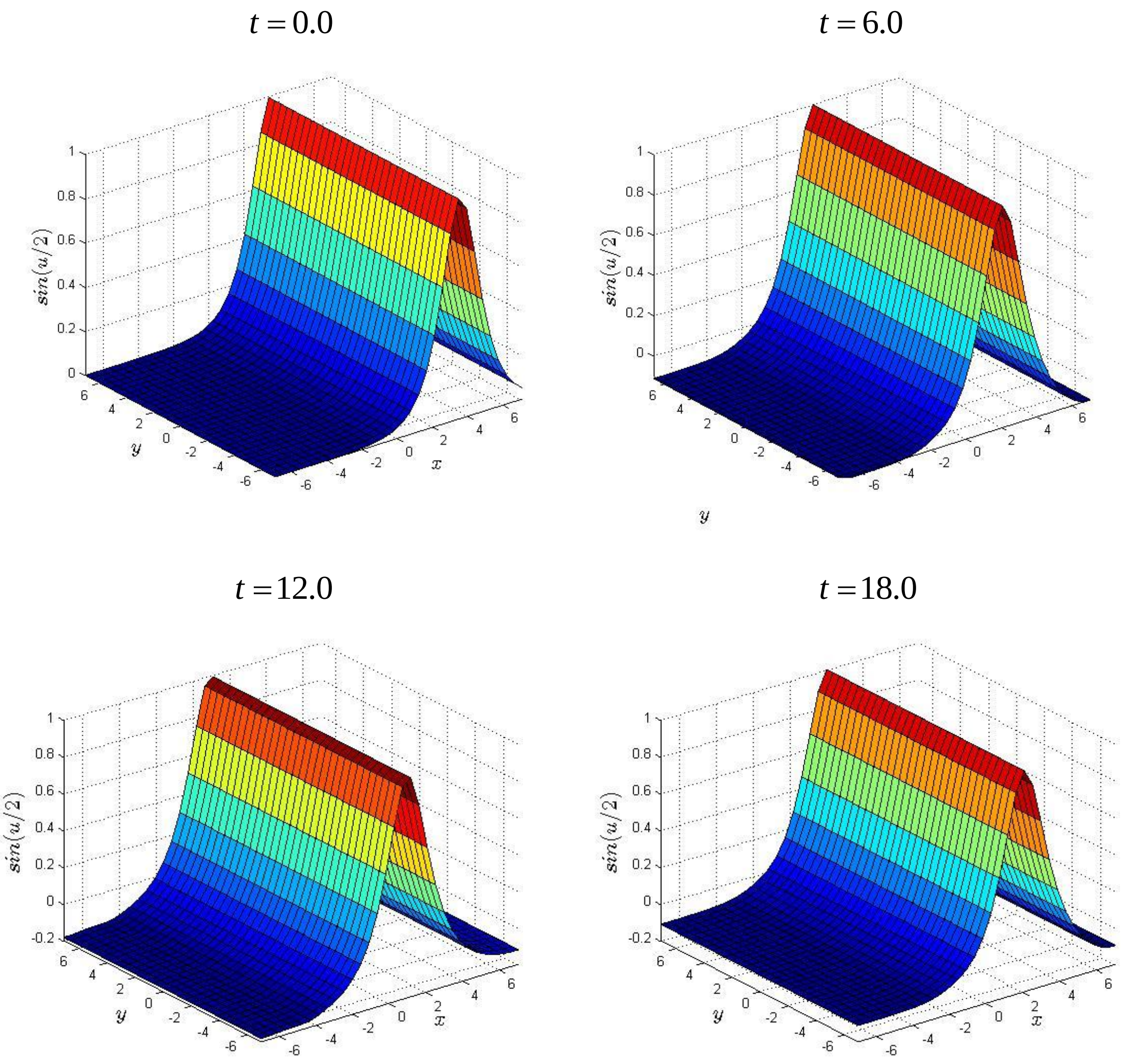

**Fig. 6.** Initial conditions and Numerical solutions of line soliton in an inhomogeneous medium at times $t = 6.0,\ 12.0,\ 18.0$ with $\Delta t = 0.001$ and grid size $31 \times 31$.

## References

1. J. K. Perring, T.H. Skyrme, A model unified field equation, Nucl. Phys. 31, (1962), 550-555.
2. G.B. Whitham, Linear and Nonlinear Waves, Wiley-Interscience, New York, NY, 1999.
3. A. Barone, F. Esposito, C. J. Magee, A.C. Scott, Theory and applications of the Sine-Gordon equation, Rivista del NuovoCimento, 1(2), (1971), 227–267.
4. J.D. Josephson, Supercurrents through barriers, Adv. Phys. 14, (1965), 419–451.
5. R.K. Dodd, I.C. Eilbeck, J.D. Gibbon, Solitons and Nonlinear Wave Equations, Academic, London, 1982.
6. B.Y. Guo, P.J. Pascual, M.J. Rodriguez, L. Vázquez, Numerical solution of the sine-Gordon equation, Appl. Math. Comput., 18, (1986), 1–14.
7. K. Djidjeli, W.G. Price, E.H. Twizell, Numerical solutions of a damped sine-Gordon equation in two space variables, J. Eng. Math., 29, (1995), 347–369.

8. P.L. Christiansen, P.S. Lomdahl, Numerical solution of (2 + 1) dimensional sine-Gordon solitons, Physica D, 2, (1981), 482–494.
9. J. Argyris, M. Haase, J.C. Heinrich, Finite element approximation to two dimensional sine-Gordon solitons, Comput. Methods Appl. Mech. Eng., 86, (1991), 1–26.
10. J. X. Xin, Modeling light bullets with the two-dimensional sine-Gordon equation, Physica D, 135, (2000), 345–368.
11. A.A. Minzoni, N. F. Smythb, A. L. Worthy, Pulse evolution for a two-dimensional sine-Gordon equation, Physica D, 159, (2001), 101–123.
12. A.A. Minzoni, N. F. Smythb, A. L. Worthy, Evolution of two-dimensional standing and travelling breather solutions for the sine–Gordon equation, Physica D, 189, (2004), 167–187.
13. Q. Sheng, A.Q.M. Khaliq, D.A. Voss, Numerical simulation of two-dimensional sine-Gordon solitons via a split cosine scheme, Math. Comput. Simulation, 68, (2005), 355–373.
14. A.G. Bratsos, An explicit numerical scheme for the sine-Gordon equation in (2+1) dimensions, Appl. Numer. Anal. Comput. Math. 2 (2), (2005), 189-211.
15. A.G. Bratsos, A modified predictor‑corrector scheme for the two dimensional sine-Gordon equation, Numer. Algor. 43, (2006), 295-308.
16. A.G. Bratsos, The solution of the two-dimensional sine-Gordon equation using the method of lines, J. Comput. Appl. Math. 206, (2007), 251-277.
17. A.G. Bratsos, A third order numerical scheme for the two-dimensional sine-Gordon equation, Mathematics and Computers in Simulation, 76, (2007), 271-282.
18. A.G. Bratsos, An improved numerical scheme for the sine-Gordon equation in (2 +1) dimensions, Int. J. Numer. Meth. Eng. 75 (2008) 787–799.
19. M. Dehghan, D. Mirzaei, The dual reciprocity boundary element method (DRBEM) for two-dimensional sine-Gordon equation, Comput. Methods Appl. Mech. Engrg., 197, (2008), 476–486.
20. M. Dehghan, A. Shokri, A numerical method for solution of the two dimensional sine-Gordon equation using the radial basis functions, Math. Comput. Simulation, 79, (2008), 700–715.
21. D. Mirzaei, M. Dehghan, Boundary element solution of the two dimensional sine-Gordon equation using continuous linear elements, Eng. Anal. Boundary Elements, 33 (2009) 12–24.
22. D. Mirzaei, M. Dehghan, Implementation of meshless LBIE method to the 2D non-linear SG problem, Int. J. Numer. Methods Engrg., 79, (2009), 1662–1682.
23. M. Dehghan, A. Ghesmati, Numerical simulation of two-dimensional sine-Gordon solitons via a local weak meshless technique based on the radial point interpolation method (RPIM), Comput. Phys. Comm., 181, (2010), 772–786.
24. D. Mirzaei, M. Dehghan, Meshless local Petrov–Galerkin (MLPG) approximation to the two-dimensional sine-Gordon equation, J. Comput. Appl. Math., 233, (2010), 2737–2754.
25. R. Jiwari, S. Pandit, R.C. Mittal, Numerical simulation of two-dimensional sine-Gordon solitons by differential quadrature method, Computer Physics Communications, 183, (2012), 600-616.
26. R. Bellman, B. G. Kashef, J. Casti, Differential quadrature: a technique for the rapid solution of nonlinear differential equations, J. Comput. Phy. 10, (1972), 40-52.
27. C. Shu, B. E. Richards, Application of generalized differential quadrature to solve two dimensional incompressible navier-Stokes equations, Int. J. Numer. Meth. Fluids, 15, (1992), 791-798.

28. J. R. Quan, C.T. Chang, New insights in solving distributed system equations by the quadrature methods-I, Comput. Chem. Eng. 13, (1989), 779–788.
29. J. R. Quan, C.T. Chang, New insights in solving distributed system equations by the quadrature methods-II, Comput. Chem. Eng. 13, (1989), 1017–1024.
30. C. Shu, Y.T. Chew, Fourier expansion-based differential quadrature and its application to Helmholtz eigenvalue problems, Commun. Numer. MethodsEng. 13 (8), (1997), 643–653.
31. C. Shu, H. Xue, Explicit computation of weighting coefficients in the harmonic differential quadrature, J. Sound Vib. 204 (3), (1997), 549–555.
32. C. Shu, Differential Quadrature and its Application in Engineering, Athenaeum Press Ltd., Great Britain, 2000)
33. A. Korkmaz, I. Dag, Cubic B-spline differential quadrature methods and stability for Burgers' equation, Eng. Comput. Int. J. Comput. Aided Eng. Software, 30 (3), (2013), 320–344.
34. A. Korkmaz, I. Dag, Shock wave simulations using sinc differential quadrature method, Eng. Comput. Int. J. Comput. Aided Eng. Software, 28(6), (2011), 654–674.
35. G. Arora, B. K. Singh, Numerical solution of Burgers' equation with modified cubic B-spline differential quadrature method, Applied Math. Comput., 224 (1), (2013), 166-177.
36. R. Jiwari, J. Yuan, A computational modeling of two dimensional reaction–diffusion Brusselator system arising in chemical processes, J. Math. Chem., J. Math. Chem., 52 (6), 1535-1551, (2014).
37. S. Gottlieb,C. W. Shu, E. Tadmor, Strong Stability-Preserving High-Order Time Discretization Methods, SIAM REVIEW, 43(1), (2001), 89-112.
38. J. R. Spiteri, S. J. Ruuth, A new class of optimal high-order strong stability-preserving time-stepping schemes, SIAM J. Numer. Analysis. 40 (2), (2002), 469-491.
39. S. Gottlieb, D. I. Ketcheson, C. W. Shu, High Order Strong Stability Preserving Time Discretizations, J. Sci. Comput.,38, (2009), 251-289.